\documentclass{article}
\usepackage{graphicx}
\usepackage{amsmath}
\usepackage{amsfonts}
\usepackage{amssymb}
\newtheorem{theorem}{Theorem}[section]

\newtheorem{corollary}[theorem]{Corollary}

\newtheorem{example}[theorem]{Example}

\newtheorem{lemma}[theorem]{Lemma}

\newtheorem{proposition}[theorem]{Proposition}
\newtheorem{remark}[theorem]{Remark}

\newenvironment{proof}[1][Proof]{\textbf{#1.} }{\ \rule{0.5em}{0.5em}}

\begin{document}

\begin{center}
\footnotetext{2000 MS Classification: 14H50.
\par
This paper is partially supported by KBN Grant 2 P03A 007 18.}{\Huge The join
of algebraic curves}

\bigskip

{\Large Tadeusz Krasi\'{n}ski}

\bigskip

28.04.2000
\end{center}

\textbf{Abstract. }{\small An effective description of the join of algebraic
curves in the complex projective space }$\mathbb{P}^{n}$ is given.

\section{Introduction}

Let $\mathbb{P}^{n}$ be the $n$-dimensional projective space over $\mathbb{C}
$. Denote by $G(1,\mathbb{P}^{n})$ the grassmannian of the all projective
lines in $\mathbb{P}^{n}.$ By the Pl\"{u}cker embedding $G(1,\mathbb{P}%
^{n})\hookrightarrow\mathbb{P}^{\binom{n+1}{2}-1}$ the grassmannian is an
algebraic subset of $\mathbb{P}^{\binom{n+1}{2}-1}.$ For any projective line
$L\subset\mathbb{P}^{n}$ we will denote by $[L]$ the corresponding point of
$G(1,\mathbb{P}^{n})$ and for any $P,Q\in\mathbb{P}^{n},$ $P\neq Q,$ we will
denote by $\overline{PQ}$ the unique projective line in $\mathbb{P}^{n}$
\ spanned by $P$ and $Q$. Likewise, for any projective subspaces
$L,K\subset\mathbb{P}^{n}$ we will denote by $\operatorname*{Span}(L,K)$ the
unique projective subspace in $\mathbb{P}^{n}$ spanned by $L$ and $K.$

If $X$ is an algebraic subset of $\mathbb{P}^{n}$ then $\operatorname*{Sing}%
(X)$ is the set of singular points of $X$. For $P\in X-\operatorname*{Sing}%
(X)$ by $T_{P}X\subset\mathbb{P}^{n}$ we denote the embedded tangent space to
$X$ at $P$.

Let $X,Y\subset\mathbb{P}^{n}$ be two varieties in $\mathbb{P}^{n}$ i.e.
irreducible algebraic subsets of $\mathbb{P}^{n}.$ The definition of the join
of $X$ and $Y$ is as follows (see \cite{H}, p.88, \cite{Z}, p.15, \cite{FOV},
Def. 1.3.5). Define the subsets of the grassmannian
\begin{align*}
\mathcal{J}^{0}(X,Y)  & :=\{[\overline{PQ}]\in G(1,\mathbb{P}^{n}):P\in X,Q\in
Y,P\neq Q\},\\
\mathcal{J}(X,Y)  & :=\overline{\mathcal{J}^{0}(X,Y)}\text{ - the closure of
}\mathcal{J}^{0}(X,Y)\text{ in }G(1,\mathbb{P}^{n})
\end{align*}
and the corresponding subsets of the projective space
\begin{align*}
J^{0}(X,Y)  & :=\bigcup_{[L]\in\mathcal{J}^{0}(X,Y)}L,\\
J(X,Y)  & :=\bigcup_{[L]\in\mathcal{J}(X,Y)}L.
\end{align*}
$\mathcal{J}(X,Y)$ and $J(X,Y)$ are algebraic subsets of $G(1,\mathbb{P}^{n})
$ and $\mathbb{P}^{n},$ respectively. $\mathcal{J}(X,Y)$ is called \textit{the
variety of lines joining \ }$X$ \textit{and} $Y$, and $J(X,Y)$ - \textit{the
join of }$X$ \textit{and} $Y$. In the case $X=Y$ the set $J(X,Y)$ is called
\textit{the secant variety of }$X$ and is denoted by $\operatorname*{Sec}(X) $
or $X^{2}$.

If $X\cap Y=\emptyset$ then we have $\mathcal{J}(X,Y)=\mathcal{J}^{0}(X,Y)$.
In the case $X\cap Y\neq\emptyset,$ the inclusion $\mathcal{J}^{0}%
(X,Y)\subset\mathcal{J}(X,Y)$ is, in general, strict. Harris in \cite{H} posed
the question which additional projective lines besides those containing points
$P\in X,Q\in Y,P\neq Q,$ are in $\mathcal{J}(X,Y)$? In the paper we give a
complete solution of this problem in the case $X,Y$ are arbitrary projective
curves (in particular for $X=Y$).

The key notion in the solution is the relative tangent cone $C_{P}(X,Y)$ to a
pair of algebraic or analytic sets $X,Y$ in a given common point $P\in X\cap
Y$ (in \cite{FOV}, S.2.5, it is denoted by $LJoin_{P}(X,Y)$). It is a
generalization of one of the Whitney's cones, precisely $C_{5}(V,P)$
(\cite{W1}, p.212, \cite{W3}, p.211), to the case of a pair of sets. The cone
$C_{P}(X,Y)$ was introduced by Achilles, Tworzewski and Winiarski \cite{ATW}
in the analytic case when $X$ and $Y$ meet at a point. This notion was used in
the new improper intersection theory in algebraic and analytic geometry (
\cite{FOV}, \cite{T}, \cite{CKT}, \cite{Cy}). It is easy to show (Proposition
\ref{jo}) that for varieties $X,Y\subset\mathbb{P}^{n}$
\[
J(X,Y)=J^{0}(X,Y)\cup\bigcup_{P\in X\cap Y}C_{P}(X,Y).
\]
So, the question is reduced to the problem of describing of $C_{P}(X,Y)$. If
$P$ is an isolated point of intersection of two analytic curves $X$ and $Y$
Ciesielska in \cite{C} proved that the cone $C_{P}(X,Y)$ is a finite sum of
two-dimensional hyperplanes. The main result of the paper (Theorem \ref{opis})
is an effective formula for the relative tangent cone $C_{P}(X,Y)$ in the
general case $X,Y$ are arbitrary analytic curves and $P\in X\cap Y$ (even in
the case $X=Y$). This formula is expressed in terms of local parametrizations
of $X$ and $Y$ at $P$. The existence of local parametrizations is the reason
for which we lead considerations over $\mathbb{C}.$

In the last section we summarize all results in Theorem \ref{last} which gives
a detailed description of the join of algebraic curves.

\section{Relative tangent cones to analytic sets}

Since the relative tangent cone is a local notion we will lead considerations
in $\mathbb{C}^{n}$ and in the case $X,Y$ are analytic sets. First we consider
the case when the point $P$ is the origin i.e. $P=\mathbf{0}$. We start from
the notion of the ordinary tangent cone to an analytic set.

Let $X$ be an analytic set in a neighbourhood $U$ of $\mathbf{0}\in
\mathbb{C}^{n}$ such that $\mathbf{0}\in X$. \textit{The tangent cone}
$C_{0}(X)$ \textit{of }$X$ \textit{at} $0$ is defined to be the set of
$\mathbf{v}\in\mathbb{C}^{n}$ with the property: there exist sequences
$(\mathbf{x}_{\nu})_{\nu\in\mathbb{N}}$ of points of $X$ and $(\lambda_{\nu
})_{\nu\in\mathbb{N}}$ of complex numbers such that
\[
\mathbf{x}_{\nu}\rightarrow0\text{ \ and \ }\lambda_{\nu}\mathbf{x}_{\nu
}\rightarrow\mathbf{v}\text{ \ when }\nu\rightarrow\infty.
\]
One can find properties of the tangent cones to analytic sets in \cite{W2},
\cite{W3}, \cite{Ch}. The tangent cone is an algebraic cone in $\mathbb{C}%
^{n}$ of dimension $\dim_{0}X.$

Let $X,Y$ be analytic subsets of a neighbourhood $U$ of $\mathbf{0}%
\in\mathbb{C}^{n}$ such that $\mathbf{0}\in X\cap Y$. \textit{The relative
tangent cone }$C_{0}(X,Y)$\textit{\ of }$X$\textit{\ and }$Y$ \textit{at}
$\mathbf{0}$ is defined to be the set of $\mathbf{v}\in\mathbb{C}^{n}$ with
the property: there exist sequences $(\mathbf{x}_{\nu})_{\nu\in\mathbb{N}} $
of points of $X,$ $(\mathbf{y}_{\nu})_{\nu\in\mathbb{N}}$ of points of $Y$ and
$(\lambda_{\nu})_{\nu\in\mathbb{N}}$ of complex numbers such that
\[
\mathbf{x}_{\nu}\rightarrow0,\text{ \ }\mathbf{y}_{\nu}\rightarrow0,\text{
\ }\lambda_{\nu}(\mathbf{y}_{\nu}-\mathbf{x}_{\nu})\rightarrow\mathbf{v}%
,\text{ \ when }\nu\rightarrow\infty.
\]

Immediately from the definition we obtain:

\begin{enumerate}
\item $C_{0}(X,Y)$ is a cone with vertex at $\mathbf{0}$.

\item  If $Y=\{\mathbf{0}\}$, then $C_{0}(X,Y)=C_{0}(X),\label{1}$

\item $C_{0}(X,Y)=C_{0}(Y,X),\label{2}$

\item $C_{0}(X,Y)$ depends only on the germs of $X$ and $Y$ at $0,\label{3}$

\item $C_{0}(X_{1}\cup X_{2},Y)=C_{0}(X_{1},Y)\cup C_{0}(X_{2},Y)$ if
$X_{1},X_{2}$ are analytic sets containing $\mathbf{0}.\label{5}$
\end{enumerate}

Next two propositions are known. Since, in the sequel, we will use facts from
the proofs we give simple and elementary proofs of them in the analytic case.
We will assume in the sequel of this section that $X,Y$ are analytic subsets
of a neighbourhood $U$ of $\mathbf{0}\in\mathbb{C}^{n}$ such that
$\mathbf{0}\in X\cap Y$.

\begin{proposition}
(\cite{ATW}, Property 2.9, in the case $X\cap Y=\{0\}$). $C_{0}(X,Y)$ is an
algebraic cone in $\mathbb{C}^{n}$.\label{alg}
\end{proposition}

\begin{proof}
By the Chow theorem it suffices to prove that $C_{0}(X,Y)$ is an analytic
subset of $\mathbb{C}^{n}$. We will apply the elementary Whitney method (
\cite{W1}, Th. 5.1, used there in the case $X=Y),$ altghough one can also use
the method of blowing-ups. Define the holomorphic functions
\begin{align*}
\alpha_{jk}  & :\mathbb{C}^{n}\times\mathbb{C}^{n}\times\mathbb{C}%
^{n}\rightarrow\mathbb{C},\mathbb{\;\;\;\;}j,k=1,...,n,\\
\alpha_{jk}(\mathbf{x,y,v})  & :=\left|
\begin{array}
[c]{cc}%
y_{j}-x_{j} & y_{k}-x_{k}\\
v_{j} & v_{k}%
\end{array}
\right|  ,
\end{align*}
where $\mathbf{x}=(x_{1},...,x_{n}),$ $\mathbf{y}=(y_{1},...,y_{n})$\ and
$\mathbf{v}=(v_{1},...,v_{n})$.

The all functions $\alpha_{jk}$ vanish if and only if $\mathbf{x}=\mathbf{y}$
or $\mathbf{v}$ is a multiple of $\mathbf{y}-\mathbf{x}$. Set
\[
B:=\{(\mathbf{x,y,v}):\mathbf{x,y\in} U,\;\alpha_{jk}(\mathbf{x,y,v}%
)=0,\;j,k=1,...,n\}.
\]
This is an analytic subset of $U\times U\times\mathbb{C}^{n}$ and hence so is
\[
B^{\prime}:=B\cap(X\times Y\times\mathbb{C}^{n}).
\]
The set $\Delta:=\{(\mathbf{x,x}):\mathbf{x}\in X\cap Y\}\subset U\times U$ is
also analytic. So,
\[
B^{\prime\prime}:=\overline{(B^{\prime}-(\Delta\times\mathbb{C}^{n}))}%
\cap(U\times U\times\mathbb{C}^{n})
\]
is an analytic set in $U\times U\times\mathbb{C}^{n}$. Then
\[
C_{0}^{\prime}(X,Y):=B^{\prime\prime}\cap(\{(\mathbf{0},\mathbf{0}%
)\}\times\mathbb{C}^{n})
\]
is analytic in $U\times U\times\mathbb{C}^{n}$. Since $\mathbf{v}\in
C_{0}(X,Y)$ if and only if $(\mathbf{0},\mathbf{0},\mathbf{v)}\in
C_{0}^{\prime}(X,Y)$, then $C_{0}(X,Y)$ is an analytic subset of
$\mathbb{C}^{n}$.
\end{proof}

\begin{proposition}
(cf. \cite{FOV}, Prop. 2.5.5) \label{dim}$\dim C_{0}(X,Y)\leqslant\dim
_{0}X+\dim_{0}Y.$
\end{proposition}

\begin{proof}
Since $C_{0}(X,Y)$ depends only on the germs of $X$ and $Y$ at $\mathbf{0}$,
we may assume that $\dim X=\dim_{0}X$ and $\dim Y=\dim_{0}Y$. Consider\ the
analytic set $B^{\prime\prime}\subset U\times U\times\mathbb{C}^{n}$, defined
in the proof of the previous Proposition. If we denote by $\pi$ the projection
$U\times U\times\mathbb{C}^{n}\rightarrow U\times U$, then $\pi(B^{\prime
\prime})\subset X\times Y$ and over each point $(\mathbf{x},\mathbf{y}%
)\in\left(  X\times Y\right)  -\Delta$ we have $(\pi|B^{\prime\prime}%
)^{-1}(\mathbf{x},\mathbf{y})=\{(\mathbf{x},\mathbf{y,}\lambda(\mathbf{y-x}%
)):\lambda\in\mathbb{C}\}$ and hence $\dim(\pi|B^{\prime\prime})^{-1}%
(\mathbf{x},\mathbf{y})=1.$ Since
\begin{equation}
B^{\prime\prime}=\overline{(\pi|B^{\prime\prime})^{-1}(X\times Y-\Delta
)},\label{w1}%
\end{equation}
then
\[
\dim B^{\prime\prime}=\dim X+\dim Y+1.
\]
By the same equality (\ref{w1}) no irreducible component of $B^{\prime\prime}$
is contained in $\Delta\times\mathbb{C}^{n}$ and in particular in
$(\mathbf{0},\mathbf{0})\times\mathbb{C}^{n}$. Hence
\[
\dim C_{0}^{\prime}(X,Y)=\dim(B^{\prime\prime}\cap(\{(\mathbf{0}%
,\mathbf{0})\}\times\mathbb{C}^{n}))\leqslant\dim X+\dim Y.
\]
\end{proof}

\begin{remark}
If we do some additional assumptions on $X$ and $Y$ then the above inequality
becomes an equality. Namely in \cite{ATW} there was proved that if $X\cap
Y=\{\mathbf{0\}}$ then $\dim C_{0}(X,Y)=\dim_{0}X+\dim_{0}Y.$ Of course, it is
no longer true in the general case.
\end{remark}

Before the next proposition we precise some notions concerning analytic
curves. By \textit{an analytic curve} we mean an analytic set $\Gamma$ of pure
dimension 1 in an open set $U\subset\mathbb{C}^{n}$. For $P\in\Gamma$ we
denote by $(\Gamma)_{P}$ the germ of $\Gamma$ at $P$ and by $\deg_{P}\Gamma$ -
\textit{the degree of }$\Gamma$ \textit{at} $P.$ \textit{A parametrization of}
$\Gamma$ \textit{at} $P$ is a holomorphic homeomorphism $\Phi:K(r)\rightarrow
U$ $(K(r):=\{z\in\mathbb{C}:|z|<r\}$ is an open disc$)$ such that $\Phi(0)=P$
and $\Phi(K(r))=\Gamma\cap U^{\prime}$ ($U^{\prime}\subset U$ is an open
neighbourhood of $P$). Then any superposition $\Phi(t^{k})$, $k\in\mathbb{N}$
we will call \textit{a description of} $X$ \textit{at }$P$. It is known that
any analytic curve $\Gamma$ such that $(\Gamma)_{P}$ is irreducible has a
parametrization. If $0\neq\Phi=(\varphi_{1},...,\varphi_{n}),$ $\Phi
(0)=\mathbf{0}$, then we define
\[
\operatorname*{ord}\Phi:=\min(\operatorname*{ord}\varphi_{1}%
,...,\operatorname*{ord}\varphi_{n}).
\]
If $\Phi$ is a parametrization of $\Gamma$ at $\mathbf{0}$ then we have
\[
\deg_{0}\Gamma=\operatorname*{ord}\Phi.
\]
It is well known that if $\Gamma$ is an analytic curve in a neighbourhood $U $
of $\mathbf{0\in}\mathbb{C}^{n}$ and $\Phi$ is its parametrization at
$\mathbf{0}$ then $C_{0}(\Gamma)$ is a line $\mathbb{C}\mathbf{v}$, where
\[
\mathbf{v}=\lim_{t\rightarrow0}\frac{\Phi(t)}{t^{\operatorname*{ord}\Phi}}.
\]
We will shortly denote this fact by
\[
\Phi(t)\underset{t\rightarrow0}{\rightsquigarrow}\mathbf{v.}%
\]
or in more condensed form $\Phi(t)\rightsquigarrow\mathbf{v}$. Note that for
any vector $\mathbf{w\in}\mathbb{C}\mathbf{v}$, by a slight change of
parameter $t\rightarrow\alpha t$, $\alpha\in\mathbb{C}$, we get that
$\Phi(\alpha t)\rightsquigarrow\mathbf{w}$. So, $\Phi$ gives rather the whole
line $\mathbb{C}\mathbf{v}$ than the vector $\mathbf{v}$ alone. So, we will
also use the notation $\Phi(t)\rightsquigarrow\mathbf{w}$ for any
$\mathbf{w\in}\mathbb{C}\mathbf{v}$.

\begin{proposition}
\label{sel}Assume that $\dim_{0}(X\cup Y)>0$. For any $0\neq\mathbf{v\in}%
C_{0}(X,Y)$ there exists an analytic curve $\Gamma\subset X\times Y$ having a
parametrization $\Phi=(\Phi_{X},\Phi_{Y}):K(r)\rightarrow X\times Y$ at
$(\mathbf{0,0)}$ such that
\[
\Phi_{Y}(t)-\Phi_{X}(t)\rightsquigarrow\mathbf{v.}%
\]
\end{proposition}

\begin{proof}
Consider the analytic set $B^{\prime\prime}\subset U\times U\times
\mathbb{C}^{n}$ defined in the proof of Proposition \ref{alg}. We have
$P:=(\mathbf{0,0,v)}\in B^{\prime\prime}.$ Since this point lies in the
closure of $B^{\prime}-(\Delta\times\mathbb{C}^{n})$ then there exists an
analytic curve $\Gamma^{\prime}\subset B^{\prime\prime}$ passing through $P$
such that $\Gamma^{\prime}-\{P\}\subset B^{\prime}-(\Delta\times\mathbb{C}%
^{n})$. Take a parametrization $(\Phi_{X}(t),\Phi_{Y}(t),\mathbf{v(}%
t\mathbf{)),}$ $t\in K(r),$ at $P$ of one irreducible component of
$(\Gamma^{\prime})_{P}.$ We have $(\Phi_{X}(0),\Phi_{Y}(0),\mathbf{v(}%
0\mathbf{))=(0,0,v).}$ Since for any $t\in K(r)$, $\Phi_{Y}(t)-\Phi_{X}(t)$
and $\mathbf{v(}t\mathbf{)}$ are linearly dependent and $\mathbf{v(}%
t\mathbf{)\rightarrow v}$ when $t\rightarrow0$ then $\Phi_{Y}(t)-\Phi
_{X}(t)\rightsquigarrow\mathbf{v.}$
\end{proof}

\begin{proposition}
(\cite{ATW}, Prop. 2.10 in the case $X\cap Y=\{\mathbf{0}\}$). $C_{0}%
(X)+C_{0}(Y)\subset C_{0}(X,Y).$
\end{proposition}

\begin{proof}
Let $0\neq\mathbf{v\in} C_{0}(X)$, $0\neq\mathbf{w\in} C_{0}(Y)$. Since
$C_{0}(X)$ is a cone then $-\mathbf{v\in} C_{0}(X)$. Take analytic curves
$\Gamma\subset X$ and $\Gamma^{\prime}\subset Y$ having parametrizations
$\Phi(t)$ and $\Psi(t)$ at $\mathbf{0}$, $t\in K(r),$ such that $\Phi
(t)\rightsquigarrow-\mathbf{v}$ and $\Psi(t)\rightsquigarrow\mathbf{w}$. Since
$\Phi(t^{\operatorname*{ord}\Psi})\in X$ and $\Psi(t^{\operatorname*{ord}\Phi
})\in Y$ for sufficiently small $t$ and
\[
\Psi(t^{\operatorname*{ord}\Phi})-\Phi(t^{\operatorname*{ord}\Psi
})\rightsquigarrow\mathbf{v+w}%
\]
then $\mathbf{v+w\in} C_{0}(X,Y)$.
\end{proof}

We will need in the sequel a propositon which was proved in \cite{ATW}, Prop.
2.10. For completness of the paper we shall give another proof of it following
easily from Proposition \ref{sel}.

\begin{proposition}
\label{trans}If $C_{0}(X)\cap C_{0}(Y)=\{\mathbf{0\}}$ then
\[
C_{0}(X,Y)=C_{0}(X)+C_{0}(Y).
\]
\end{proposition}

\begin{proof}
It suffices to prove
\[
C_{0}(X,Y)\subset C_{0}(X)+C_{0}(Y).
\]
Take $0\neq\mathbf{w\in} C_{0}(X,Y)$. We may assume that $\mathbf{w\notin}%
C_{0}(X)\cup C_{0}(Y)$. By Proposition \ref{sel} there exists an analytic
curve $\Gamma\subset X\times Y$ having a parametrization $\Phi=(\Phi_{X}%
,\Phi_{Y}):K(r)\rightarrow X\times Y$ at $(\mathbf{0,0)}$ such that
\[
\Phi_{Y}(t)-\Phi_{X}(t)\rightsquigarrow\mathbf{w.}%
\]
Since $\mathbf{w\notin} C_{0}(X)$ and $\mathbf{w\notin} C_{0}(Y)$ then
\begin{equation}
\operatorname*{ord}\Phi_{Y}=\operatorname*{ord}\Phi_{X}<+\infty.\label{ddd}%
\end{equation}
Let
\begin{align*}
\Phi_{X}(t)  & \rightsquigarrow\mathbf{v}_{1},\;\;\;0\neq\mathbf{v}%
_{1}\mathbf{\in} C_{0}(X),\\
\Phi_{Y}(t)  & \rightsquigarrow\mathbf{v}_{2},\;\;\;0\neq\mathbf{v}%
_{2}\mathbf{\in} C_{0}(Y).
\end{align*}
Since $C_{0}(X)\cap C_{0}(Y)=\{\mathbf{0\}}$ then $\mathbf{v}_{1}$and
$\mathbf{v}_{2}$ are linearly independent. Hence and from (\ref{ddd})
\[
\Phi_{Y}(t)-\Phi_{X}(t)\rightsquigarrow\mathbf{v}_{2}-\mathbf{v}_{1}\mathbf{.}%
\]
So, $\mathbf{w=v}_{2}-\mathbf{v}_{1}\in C_{0}(X)+C_{0}(Y).$
\end{proof}

Let now $X,Y$ be analytic subsets of a neighbourhood $U$ of a point
$P\in\mathbb{C}^{n}$ such that $P\in X\cap Y$. We define the \textit{relative
tangent cone }$C_{P}(X,Y)$\textit{\ of }$X$\textit{\ and }$Y$\textit{\ at }$P
$ by
\[
C_{P}(X,Y):=P+C_{0}(X-P,Y-P)
\]

\section{Relative tangent cone to analytic curves}

In the case $X,Y$ are analytic curves we may give a more detailed description
of $C_{0}(X,Y).$ The aim of this section is to give an effective formula for
$C_{0}(X,Y)$ in terms of local parametrizations of $X$ and $Y.$

First, we formulate a useful lemma which is a a simple generalization of
Proposition \ref{sel}.

\begin{lemma}
\label{cur}Let $X,Y$ be analytic curves in a neighbourhood of $\mathbf{0\in
}\mathbb{C}^{n}$ such that $\mathbf{0}\in X\cap Y$ and the germs
$(X)_{\mathbf{0}},(Y)_{\mathbf{0}}$ are irreducible. Let $\Phi(t)$ and
$\Psi(\tau),$ $t,\tau\in K(r),$ be parametrizations of $X$ and $Y$ at
$\mathbf{0}$. Then for any $\mathbf{v\in} C_{0}(X,Y)$ there exists an analytic
curve $\Gamma\subset K(r)\times K(r)$ having a parametrization $\Theta
(s)=(t(s),\tau(s)):K(r^{\prime})\rightarrow K(r)\times K(r)$ at
$(\mathbf{0,0)}$ such that
\[
\Phi(t(s))-\Psi(\tau(s))\rightsquigarrow\mathbf{v.}%
\]
Moreover, we have the same result if $\Phi$ and $\Psi$ are only descriptions
of $X$ and $Y$ at $\mathbf{0}$.
\end{lemma}

\begin{proof}
The proof follows from Proposition \ref{sel} and the fact that the mapping
$(\Phi$,$\Psi)$ is an analytic cover.
\end{proof}

Now we prove a key proposition for a description of relative tangent cones.
This proposition was proved by Ciesielska \cite{C} in the case $X\cap
Y=\{\mathbf{0\}}$, although the idea of her proof can be used in the more
general case $\mathbf{0}\in X\cap Y$.

\begin{proposition}
\label{ert}Let $X,Y$ be analytic curves in a neighbourhood of $\mathbf{0\in
}\mathbb{C}^{n}$ such that $\mathbf{0}\in X\cap Y$. Then
\[
C_{0}(X,Y)+C_{0}(X)=C_{0}(X,Y).
\]
\end{proposition}

\begin{proof}
We may assume that the germs $(X)_{\mathbf{0}},(Y)_{\mathbf{0}}$ are
irreducible. It suffices to prove that
\begin{equation}
C_{0}(X,Y)+C_{0}(X)\subset C_{0}(X,Y).\label{12}%
\end{equation}
Since $X,Y$ are analytic curves and $(X)_{\mathbf{0}},(Y)_{\mathbf{0}}$ are
irreducible at $\mathbf{0}$ we will consider two possible cases:

1$^{\circ}.$ $C_{0}(X)\cap C_{0}(Y)=\{\mathbf{0}\}$. Then by Proposition
\ref{trans} $C_{0}(X,Y)=C_{0}(X)+C_{0}(Y)$. Hence we get (\ref{12}).

2$^{\circ}.$ $C_{0}(X)=C_{0}(Y)$. After a linear change of coordinates in
$\mathbb{C}^{n}$ we may assume that $C_{0}(X)=\mathbb{C}\mathbf{e}_{1}$, where
$\mathbf{e}_{1}:=(1,0,...,0).$ Put $k:=\deg_{0}X,$ $l:=\deg_{0}Y.$ Let $\Phi$
and $\Psi$ be parametrizations of $X$ and $Y$ at $\mathbf{0}$,
respectively$\mathbf{.}$ Since $C_{0}(X)=C_{0}(Y)=\mathbb{C}\mathbf{e}_{1}$,
we may assume that
\begin{align}
\Phi(t)  & =(t^{k},\phi_{2}(t),...,\phi_{n}(t)),\;t\in
K(r),\;\operatorname*{ord}\phi_{i}>k,\;\;\;i=2,...,n,\label{fis}\\
\Psi(\tau)  & =(\tau^{l},\psi_{2}(\tau),...,\psi_{n}(\tau)),\;\tau\in
K(r),\;\operatorname*{ord}\psi_{i}>l,\;\;\;i=2,...,n.\label{psis}%
\end{align}
Consider descriptions of $X$ and $Y$
\begin{align*}
\tilde{\Phi}(t)  & :=\Phi(t^{l})=(t^{kl},\phi_{2}(t^{l}),...,\phi_{n}%
(t^{l})),\;\;t\in K(\tilde{r}),\\
\tilde{\Psi}(\tau)  & :=\Psi(\tau^{k})=(\tau^{kl},\psi_{2}(\tau^{k}%
),...,\psi_{n}(\tau^{k})),\;\;\tau\in K(\tilde{r}),
\end{align*}
where $\tilde{r}$ is a sufficiently small positive number.

Take now $\mathbf{0}\neq\mathbf{v}=(v_{1},...,v_{n})\mathbf{\in} C_{0}(X,Y)$
and $\mathbf{w=(}w\mathbf{,}0,...,0)\mathbf{\in} C_{0}(X)$. From Lemma
\ref{cur} there is an analytic curve $\Gamma\subset K(\tilde{r})\times
K(\tilde{r})$ having a parametrization $\Theta(s)=(t(s),\tau(s)):K(r^{\prime
})\rightarrow K(\tilde{r})\times K(\tilde{r})$ at $(\mathbf{0,0)}$ such that
\[
\tilde{\Phi}(t(s))-\tilde{\Psi}(\tau(s))\rightsquigarrow\mathbf{v.}%
\]
Define
\[
N:=\operatorname*{ord}(\tilde{\Phi}(t(s))-\tilde{\Psi}(\tau(s))).
\]
Then
\[
\mathbf{v}=\lim_{s\rightarrow0}\frac{\tilde{\Phi}(t(s))-\tilde{\Psi}(\tau
(s))}{s^{N}}.
\]
Since $\Theta$ is a parametrization of a curve we have that $t(s)$ or
$\tau(s)$ is not identically zero. Without loss of generality, we may assume
that $t(s)\not \equiv0$ and $\operatorname*{ord}t(s)\leqslant
\operatorname*{ord}\tau(s)$. Put $p:=\operatorname*{ord}t(s)$. Hence
$N\geqslant pkl$. Changing unessentially $t(s)$ we may assume that
$t(s)=s^{p}$. We define
\[
\tilde{t}(s):=s^{p}+\frac{w}{kl}s^{p+N-pkl},
\]
We claim that
\[
\tilde{\Phi}(\tilde{t}(s))-\tilde{\Psi}(\tau(s))\rightsquigarrow
\mathbf{v}+\mathbf{w.}%
\]
In fact, for the first coordinate we have
\[
\lim_{s\rightarrow0}\frac{(\tilde{t}(s))^{kl}-(\tau(s))^{kl}}{s^{N}}%
=\lim_{s\rightarrow0}\frac{(\tilde{t}(s))^{kl}-(t(s))^{kl}+(t(s))^{kl}%
-(\tau(s))^{kl}}{s^{N}}=w+v_{1}%
\]
and for the next coordinates
\begin{align*}
\lim_{s\rightarrow0}\frac{(\phi_{i}(\tilde{t}(s)^{l})-\psi_{i}(\tau(s)^{k}%
)}{s^{N}}  & =\\
\lim_{s\rightarrow0}\frac{\phi_{i}(\tilde{t}(s)^{l})-\phi_{i}(t(s)^{l}%
)+\phi_{i}(t(s)^{l})-\psi_{i}(\tau(s)^{k})}{s^{N}}  & =v_{i},\;i=2,...,n.
\end{align*}
\end{proof}

From this proposition we obtain the first description of relative tangent
cones to analytic curves (cf. \cite{C}, Cor. 3.2).

\begin{corollary}
Let $X,Y$ be analytic curves in a neighbourhoodof $\mathbf{0\in}\mathbb{C}%
^{n}$such that $\mathbf{0}\in X\cap Y$ and $(X)_{\mathbf{0}},(Y)_{\mathbf{0}}
$ be irreducible germs at $\mathbf{0}$. Then two cases may occur:

1. $C_{0}(X,Y)=C_{0}(X)=C_{0}(Y)$.

2. $C_{0}(X,Y)$ is a finite sum of two-dimensional hyperplanes.
\end{corollary}

\begin{proof}
If $C_{0}(X)\cap C_{0}(Y)=\{\mathbf{0}\}$, then by Proposition \ref{trans}
$C_{0}(X,Y)=C_{0}(X)+C_{0}(Y)$ is a two-dimensional hyperplane. If
$C_{0}(X)=C_{0}(Y),$ then taking an $(n-1)$-dimensional hyperplane $H$ through
$\mathbf{0}$, transversal to $C_{0}(X),$ we easily obtain from Proposition
\ref{ert} that
\begin{equation}
C_{0}(X,Y)=C_{0}(X,Y)\cap H+C_{0}(X).\label{hh}%
\end{equation}
Since by Proposition \ref{dim} $\dim C_{0}(X,Y)\leqslant2$ then by (\ref{hh})
$\dim C_{0}(X,Y)\cap H\leqslant1$. But $C_{0}(X,Y)\cap H$ is also an algebraic
cone. Hence $C_{0}(X,Y)\cap H$ is either $\{\mathbf{0\}}$ or a finite number
of lines. So, by (\ref{hh}), $C_{0}(X,Y)=C_{0}(X)$ in the first case or is a
finite sum of two-dimensional hyperplanes in the second one.
\end{proof}

Now we give the main result of the paper. It is a formula for the $C_{0}(X,Y)
$ in terms of parametrizations of $X$ and $Y$. First we fix some notations. By
$\mathbf{e}_{1}\mathbf{,...,e}_{n}$ we denote the versors of axes in
$\mathbb{C}^{n}$. For vectors $\mathbf{v,w}\in\mathbb{C}^{n}$ by
$\operatorname*{Lin}(\mathbf{v,w})$ we denote the hyperplane in $\mathbb{C}%
^{n}$ generated by $\mathbf{v}$ and $\mathbf{w}$\textbf{.} By
$\operatorname*{in}(\chi(s))$ of a power series $\chi(s)\not \equiv0$ we mean
its \textit{initial form} i.e. if $\chi(s)=\beta_{p}s^{p}+...,$ $\ \ \beta
_{p}\neq0,$ then $\operatorname*{in}(\chi(s))=\beta_{p}s^{p}$ (additionally we
put $\operatorname*{in}(0):=0$)$.$

\begin{theorem}
\label{opis}Let $X,Y$ be analytic curves in a neighbourhood $U$ of the point
$\mathbf{0}\in\mathbb{C}^{n}$ such that $\mathbf{0}\in X\cap Y$ and $(X)_{0}
$, $(Y)_{0}$ are irreducible germs. Let
\begin{align}
\Phi(t)  & =(t^{k},\phi_{2}(t),...,\phi_{n}(t)),\;t\in
K(r),\;\operatorname*{ord}\phi_{i}>k,\;i=2,...,n\label{fi}\\
\Psi(\tau)  & =(\tau^{l},\psi_{2}(\tau),...,\psi_{n}(\tau)),\;\tau\in
K(r),\;\operatorname*{ord}\psi_{i}>l,\;i=2,...,n\label{psi}%
\end{align}
be parametrizations of $X$ and $Y$at $\mathbf{0}.$ Assume that $l\leqslant k$.
Let $\varepsilon_{1},...,\varepsilon_{l}$ be the all roots of unity of degree
$l$. For $i=1,...,l$ we define
\begin{align*}
n_{i}  & :=\left\{
\begin{array}
[c]{ccc}%
\operatorname*{ord}(\Phi(t^{l})-\Psi(\varepsilon_{i}t^{k})) & \text{if} &
\Phi(t^{l})-\Psi(\varepsilon_{i}t^{k})\not \equiv0\\
0 & \text{if} & \Phi(t^{l})-\Psi(\varepsilon_{i}t^{k})\equiv0
\end{array}
\right.  ,\\
\mathbf{v}_{i}  & :=\lim_{t\rightarrow0}\frac{\Phi(t^{l})-\Psi(\varepsilon
_{i}t^{k})}{t^{n_{i}}}.
\end{align*}
Then
\[
C_{0}(X,Y)=\operatorname*{Lin}(\mathbf{v}_{1}\mathbf{,e}_{1})\cup
...\cup\operatorname*{Lin}(\mathbf{v}_{l}\mathbf{,e}_{1}).
\]
\end{theorem}

\begin{proof}
Instead of the parametrizations $\Phi$ and $\Psi$, we shall use descriptions
of $X$ and $Y$ . Define
\begin{align*}
\tilde{\Phi}(t)  & :=\Phi(t^{l})=(t^{kl},\phi_{2}(t^{l}),...,\phi_{n}%
(t^{l})),\;\;\;t\in K(r^{1/l}),\\
\tilde{\Psi}(\tau)  & :=\Psi(\tau^{k})=(\tau^{kl},\psi_{2}(\tau^{k}%
),...,\psi_{n}(\tau^{k})),\;\;\tau\in K(r^{1/k}).
\end{align*}
Obviously, $(\tilde{\Phi}(K(r^{1/l}))_{0}=(X)_{0}$, $(\tilde{\Psi}%
(K(r^{1/k}))_{0}=(Y)_{0}$. From the form of $\tilde{\Phi}$ and $\tilde{\Psi}$
we see that
\[
C_{0}(X)=C_{0}(Y)=\mathbb{C}\mathbf{e}_{1}.
\]
Take the hyperplane
\[
H:=\{(x_{1},...,x_{n})\in\mathbb{C}^{n}:x_{1}=0\},
\]
transversal to $C_{0}(X)=C_{0}(Y).$ From Proposition \ref{ert} we easily
obtain
\[
C_{0}(X,Y)=C_{0}(X,Y)\cap H+C_{0}(X).
\]
Since $C_{0}(X,Y)$ is an analytic cone in $\mathbb{C}^{n}$ of dimension
$\leqslant2$, then from this equality $C_{0}(X,Y)\cap H$ is either
$\{\mathbf{0}\}$ or a finite system of lines. So, it suffices to prove that
\[
C_{0}(X,Y)\cap H=\bigcup_{i=1}^{l}\mathbb{C}\mathbf{v}_{i}.
\]
By definition of $\mathbf{v}_{i}$ we have obviously
\[
\bigcup_{i=1}^{l}\mathbb{C}\mathbf{v}_{i}\subset C_{0}(X,Y)\cap H.
\]
Take now any vector $0\neq\mathbf{w}\in C_{0}(X,Y)\cap H$. By Lemma \ref{cur}
there exists an analytic curve $\Gamma\subset K(r^{1/l})\times K(r^{1/k})$
having a parametrization $\Theta(s)=(t(s),\tau(s)):K(r^{\prime})\rightarrow
K(r^{1/l})\times K(r^{1/k})$ at $(\mathbf{0,0)}$ such that
\[
\left(  \tilde{\Phi}(t(s))-\tilde{\Psi}(\tau(s))\right)  \rightsquigarrow
\mathbf{w},\text{ when }s\rightarrow0,
\]
i.e.
\[
(t(s)^{kl}-\tau(s)^{kl},\phi_{2}(t(s)^{l})-\psi_{2}(\tau(s)^{k}),...,\phi
_{n}(t(s)^{l})-\psi_{n}(\tau(s)^{k}))\rightsquigarrow\mathbf{w},\text{ when
}s\rightarrow0.
\]
Since $t(s)\not \equiv0$ or $\tau(s)\not \equiv0$ we may assume that
$t(s)\not \equiv0.$ Changing unessentially the parameter $s$ we may assume
that
\[
t(s)=s^{p},\;p\in\mathbb{N}\text{.}%
\]
Then
\[
(s^{pkl}-\tau(s)^{kl},\phi_{2}(s^{pl})-\psi_{2}(\tau(s)^{k}),...,\phi
_{n}(s^{pl})-\psi_{n}(\tau(s)^{k}))\rightsquigarrow\mathbf{w},\text{ when
}s\rightarrow0.
\]
Since $\mathbf{w}=(0,w_{2},...,w_{n})\neq0$, then there exists $j\in
\{2,...,n\}$ such that
\begin{equation}
\operatorname*{ord}(\phi_{j}(s^{pl})-\psi_{j}(\tau(s)^{k})<\operatorname*{ord}%
(s^{pkl}-\tau(s)^{kl}).\label{2}%
\end{equation}
Denote by $J$ the set of $j\in\{2,...,n\}$ for which the above inequality
holds$.$ Since $\operatorname*{ord}\phi_{j}>k$ and $\operatorname*{ord}%
\psi_{j}>l,$ then from the above inequality we obtain that $\tau(s)$ has the
form
\[
\tau(s)=\alpha_{p}s^{p}+\alpha_{p+1}s^{p+1}+...,\;\;\;\;\;\alpha_{p}^{kl}=1.
\]
Hence $\alpha_{p}^{k}=\varepsilon_{i_{0}}$ for some $i_{0}\in\{1,...,l\}.$ We
shall show that $\mathbf{w}=\mathbf{v}_{i_{0}}.$ Consider the cases:

1. the coefficients $\alpha_{r}$ vanish for $r>p$ i.e. $\tau(s)=\alpha
_{p}s^{p}.$ Then $\tau(s)^{k}=\alpha_{p}^{k}s^{pk}=\varepsilon_{i_{0}}s^{pk}.$
Hence we have $\mathbf{w}=\mathbf{v}_{i_{0}},$

2. not all the coefficients $\alpha_{r}$ vanish for $r>p.$ Let $m$ be the
smallest positive integer such that $\alpha_{p+m}\neq0.$ Then
\begin{align}
\tau(s)  & =\alpha_{p}s^{p}+\alpha_{p+m}s^{p+m}+....\nonumber\\
\tau(s)^{k}  & =\varepsilon_{i_{0}}s^{pk}+\alpha s^{pk+m}+....,\;\;\alpha
\neq0,\label{00}\\
\operatorname*{ord}(s^{pkl}-\tau(s)^{kl})  & =pkl+m,\label{33}\\
\operatorname*{ord}(\phi_{j}(s^{pl})-\psi_{j}(\tau(s)^{k})  & <pkl+m\text{
\ \ for \ \ }j\in J\label{400}\\
\operatorname*{ord}(\phi_{j}(s^{pl})-\psi_{j}(\tau(s)^{k})  & \geqslant
pkl+m\text{ \ \ for \ \ }j\notin J\label{10}%
\end{align}
Let us first note that for $j\in\{2,...,n\}$ from (\ref{00}) and the fact that
$\operatorname*{ord}\psi_{j}>l$ we have
\begin{equation}
\operatorname*{ord}(\psi_{j}(\tau(s)^{k})-\psi_{j}(\varepsilon_{i_{0}}%
s^{pk}))\geqslant pkl+m.\label{11}%
\end{equation}
Hence and from (\ref{400}) for $j\in J$ we have
\begin{align}
\operatorname*{in}\left(  \phi_{j}(s^{pl})-\psi_{j}(\tau(s)^{k}\right)   &
=\operatorname*{in}\left(  \phi_{j}(s^{pl})-\psi_{j}(\varepsilon_{i_{0}}%
s^{pk})+\psi_{j}(\varepsilon_{i_{0}}s^{pk})-\psi_{j}(\tau(s)^{k}\right)
\nonumber\\
& =\operatorname*{in}\left(  \phi_{j}(s^{pl})-\psi_{j}(\varepsilon_{i_{0}%
}s^{pk})\right)  ,\label{21}%
\end{align}
and for $j\notin J$ from (\ref{10}) we get
\begin{align}
\operatorname*{ord}\left(  \phi_{j}(s^{pl})-\psi_{j}(\varepsilon_{i_{0}}%
s^{pk})\right)   & =\nonumber\\
\operatorname*{ord}\left(  \phi_{j}(s^{pl})-\psi_{j}(\tau(s)^{k})+\psi
_{j}(\tau(s)^{k})-\psi_{j}(\varepsilon_{i_{0}}s^{pk})\right)   & \geqslant
pkl+m.\label{23}%
\end{align}
Hence
\begin{equation}
\operatorname*{ord}\left(  \Phi(s^{pl})-\Psi(\tau(s)^{k}\right)
=\operatorname*{ord}\left(  \Phi(s^{pl})-\Psi(\varepsilon_{i_{0}}%
s^{pk})\right)  =pn_{i_{0}}.\label{22}%
\end{equation}
\qquad Now, we have
\begin{align*}
\mathbf{v}_{i_{0}}  & =\lim_{t\rightarrow0}t^{-n_{i_{0}}}\left(  \Phi
(t^{l})-\Psi(\varepsilon_{i_{0}}t^{k})\right) \\
& =\lim_{s\rightarrow0}s^{-pn_{i_{0}}}\left(  \Phi(s^{pl})-\Psi(\varepsilon
_{i_{0}}s^{pk})\right) \\
& =\lim_{s\rightarrow0}s^{-pn_{i_{0}}}(0,\phi_{2}(s^{pl})-\psi_{2}%
(\varepsilon_{i_{0}}s^{pk}),...,\phi_{n}(s^{pl})-\psi_{n}(\varepsilon_{i_{0}%
}s^{pk}))\\
& =\lim_{s\rightarrow0}s^{-pn_{i_{0}}}\left(  0,\operatorname*{in}\left(
\phi_{2}(s^{pl})-\psi_{2}(\varepsilon_{i_{0}}s^{pk})\right)
,...,\operatorname*{in}\left(  \phi_{n}(s^{pl})-\psi_{n}(\varepsilon_{i_{0}%
}s^{pk})\right)  \right)  .
\end{align*}
On the other hand, from definition of $\mathbf{w}$ and (\ref{22}) we have
\begin{align*}
\mathbf{w}  & =\lim_{s\rightarrow0}\frac{\left(  \Phi(s^{pl})-\Psi(\tau
(s)^{k})\right)  }{s^{\operatorname*{ord}\left(  \Phi(s^{pl})-\Psi(\tau
(s)^{k})\right)  }}\\
& =\lim_{s\rightarrow0}s^{-pn_{i_{0}}}\left(  \Phi(s^{pl})-\Psi(\tau
(s)^{k})\right) \\
& =\lim_{s\rightarrow0}s^{-pn_{i_{0}}}(s^{pkl}-\tau(s)^{k},\phi_{2}%
(s^{pl})-\psi_{2}(\tau(s)^{k}),...,\phi_{n}(s^{pl})-\psi_{n}(\tau(s)^{k}))\\
& =\lim_{s\rightarrow0}s^{-pn_{i_{0}}}\left(  \operatorname*{in}\left(
s^{pkl}-\tau(s)^{k}\right)  ,\operatorname*{in}\left(  \phi_{2}(s^{pl}%
)-\psi_{2}(\tau(s)^{k})\right)  ,...,\operatorname*{in}\left(  \phi_{n}%
(s^{pl})-\psi_{n}(\tau(s)^{k})\right)  \right)
\end{align*}

Then from (\ref{33}), (\ref{21}), (\ref{10}), (\ref{23}) we finally obtain
\[
\mathbf{v}_{i_{0}}=\mathbf{w}.
\]

This ends the proof.
\end{proof}

\begin{remark}
From forms (\ref{fi}), (\ref{psi}) of parametrizations it follows that
$C_{0}(X)=C_{0}(Y)=\mathbb{C}\mathbf{e}_{1}.$ By Proposition \ref{trans} we
see that only this case is interesting.Moreover, the assumption on the form of
parametrizations is not restrictive, because it is well-known that for any
analytic curve $X$ with irreducible germ at $0$ there exists a linear change
of coordinates in $\mathbb{C}^{n}$ such that in the new coordinates
$C_{0}(X)=\mathbb{C}\mathbf{e}_{1}$ and there exists a parametrization of $X$
at $0$ of form (\ref{fi}).
\end{remark}

\begin{remark}
It is easily seen that by an unessential change of the parameter
$t\rightarrow\mu t$, $\mu^{kl}=1,$ for each of the vectors
\[
\mathbf{v}_{\varepsilon,\eta}:=\lim_{t\rightarrow0}\frac{\Phi(\eta t^{l}%
)-\Psi(\varepsilon t^{k})}{t^{\operatorname*{ord}(\Phi(\eta t^{l}%
)-\Psi(\varepsilon t^{k}))}},\;\;\varepsilon^{l}=1,\;\eta^{k}=1
\]
there exists $i\in\{1,...,l\}$ such that
\[
\mathbb{C}\mathbf{v}_{\varepsilon,\eta}=\mathbb{C}\mathbf{v}_{i}.
\]
\end{remark}

\begin{corollary}
Let $X,Y$ be analytic curves in a neighbourhood $U$ of the point
$\mathbf{0}\in\mathbb{C}^{n}$ such that $\mathbf{0}\in X\cap Y$ and $(X)_{0}$,
$(Y)_{0} $ are irreducible germs. Then

1. if $(X)_{0}=(Y)_{0}$ and this germ is nonsingular, then
\[
C_{0}(X,X)=C_{0}(X)=T_{P}X,
\]

2. in the remaining cases $C_{0}(X,Y)$ is the sum of $r$ two-dimensional
hyperplanes, where
\[
1\leqslant r\leqslant\min(\deg_{0}X,\deg_{0}Y)
\]
\end{corollary}

\begin{proof}
It follows from Theorem \ref{opis} by considering parametrizations of $X$ and
$Y$ at $0$ in the nonsingular case and singular one.
\end{proof}

\begin{example}
Let
\begin{align*}
X  & :=\{(t^{2},t^{3},0):t\in\mathbb{C}\}\subset\mathbb{C}^{3},\\
Y  & :=\{(\tau^{2},0,\tau^{3}):\tau\in\mathbb{C}\}\subset\mathbb{C}^{3}.
\end{align*}
$X$ and $Y$ satisfy assumptions of Theorem \ref{opis}. We have $k=l=2$ and
$\mathbf{v}_{1}=[0,1,1]$, $\mathbf{v}_{2}=[0,1,-1].$ Hence
\[
C_{0}(X,Y)=\operatorname*{Lin}(\mathbf{v}_{1},\mathbf{e}_{1})\cup
\operatorname*{Lin}(\mathbf{v}_{2},\mathbf{e}_{1})=\{(x,y,z)\in\mathbb{C}%
^{3}:y^{2}-z^{2}=0\}.
\]
\end{example}

\section{Join of algebraic curves}

In this section we answer the question posed in the introduction: which
additional projective lines besides those containing points $P\in X,Q\in
Y,P\neq Q,$ are in $\mathcal{J}(X,Y)$ in the case $X,Y$ are algebraic curves?
First, we give a relation between the join of arbitrary varieties and relative
tangent cones.

Let $X,Y$ be arbitrary algebraic subsets of $\mathbb{P}^{n}$ and $P\in X\cap
Y$. Let $U\subset\mathbb{P}^{n}$ be a canonical affine part of $\mathbb{P}%
^{n}$ such that $P\in U,$ and $\varphi:U\rightarrow\mathbb{C}^{n}$ the
corresponding canonical map. Then we define relative tangent cone $C_{P}(X,Y)
$ to $X$ and $Y$ at $P$ by
\[
C_{P}(X,Y):=\overline{\varphi^{-1}(C_{\varphi(P)}(\varphi(X\cap U),\varphi
(Y\cap U))}.
\]
One can easily check that it does not depend on the choice of the canonical
affine part $U$ of $\mathbb{P}^{n}$ (in \cite{FOV}, Def. 4.3.6, there is
another equivalent definition of $C_{P}(X,Y)$ using the affine cones $\hat
{X},\hat{Y}\subset\mathbb{C}^{n+1}$ generated by $X$ and $Y$)$.$

Since $C_{P}(X,Y)$ is a sum of projective lines passing through $P$ we may
define
\[
\mathcal{C}_{P}(X,Y):=\{[L]\in G(1,\mathbb{P}^{n}):L\subset C_{P}(X,Y)\text{
and }P\in L\}.
\]

\begin{proposition}
\label{jo}Let $X,Y$ be arbitrary algebraic subsets of $\mathbb{P}^{n}$. Then
\begin{align*}
\mathcal{J}(X,Y)  & =\mathcal{J}^{0}(X,Y)\cup\bigcup_{P\in X\cap Y}%
\mathcal{C}_{P}(X,Y),\\
J(X,Y)  & =J^{0}(X,Y)\cup\bigcup_{P\in X\cap Y}C_{P}(X,Y).
\end{align*}
\end{proposition}

\begin{proof}
Note that the topology in $G(1,\mathbb{P}^{n})$ can be described in the
following elementary way: if $[L],[L_{i}]\in G(1,\mathbb{P}^{n})$, $i=1,2,...
$, then $[L_{i}]\rightarrow\lbrack L]$ when $i\rightarrow\infty$ in
$G(1,\mathbb{P}^{n})$ if and only if there exist points $P_{i},Q_{i}\in L_{i}%
$, $i=1,2,...,$ $P_{i}\neq Q_{i}$, $P,Q\in L,$ $P\neq Q$, and their
homogeneous coordinates $P_{i}=(x_{0}^{i}:...:x_{n}^{i})$, $Q_{i}=(y_{0}%
^{i}:...:y_{n}^{i})$, $P=(x_{0}:...:x_{n})$, $Q=(y_{0}:...:y_{n}) $ such that
$x_{j}^{i}\rightarrow x_{j}$ and $y_{j}^{i}\rightarrow y_{j}$ when
$i\rightarrow\infty$ in $\mathbb{C}\ $for $j=0,1,...,n.$

Take $[L]\in\mathcal{J}(X,Y)-\mathcal{J}^{0}(X,Y)$. Then there exist
$[\overline{P_{i}Q_{i}}]\in G(1,\mathbb{P}^{n})$, $i=1,2,...$, $P_{i}\in X,$
$Q_{i}\in Y$, $P_{i}\neq Q_{i}$, such that $[\overline{P_{i}Q_{i}}%
]\rightarrow\lbrack L]$ when $i\rightarrow\infty$. Since $X,Y$ are compact
sets we may assume that $P_{i}\rightarrow P\in X$ and $Q_{i}\rightarrow Q\in
Y.$ Since $[L]\notin\mathcal{J}^{0}(X,Y)$ then $P=Q$. Hence $P\in X\cap Y$. Of
course $P\in L.$ From the above description of topology in $G(1,\mathbb{P}%
^{n})$ we easily obtain that $L\subset C_{P}(X,Y). $

The opposite inclusion $\bigcup_{P\in X\cap Y}\mathcal{C}_{P}(X,Y)\subset
\mathcal{J}(X,Y)$ is obvious.
\end{proof}

From the above proposition and the previous results we obtain the full
description of the join of algebraic curves in $\mathbb{P}^{n}.$

\begin{theorem}
\label{last}Let $X,Y$ be irreducible curves in $\mathbb{P}^{n}.$ Then:

1. if $X=Y$ then
\begin{align*}
\mathcal{J}(X,X)  & =\mathcal{J}^{0}(X,X)\cup\bigcup_{P\in\text{Sing(}%
X)}\mathcal{C}_{P}(X,X)\cup\bigcup_{P\in X-\text{Sing(}X)}[T_{P}(X)],\\
J(X,X)  & =J^{0}(X,X)\cup\bigcup_{P\in\text{Sing(}X)}C_{P}(X,X)\cup
\bigcup_{P\in X-\text{Sing(}X)}T_{P}(X),
\end{align*}

2. if $X\neq Y$ and $X\cap Y=\{P_{1},...,P_{k}\}$ then
\begin{align*}
\mathcal{J}(X,Y)  & =\mathcal{J}^{0}(X,Y)\cup\bigcup_{i=1}^{k}\mathcal{C}%
_{P_{i}}(X,Y),\\
J(X,Y)  & =J^{0}(X,Y)\cup\bigcup_{i=1}^{k}C_{P_{i}}(X,Y).
\end{align*}
\qquad Moreover, in both cases each $C_{P}(X,Y)$ is a finite sum of projective
two-dimensional hyperplanes passing through $P.$ They are effectively
described in the following way: for a given point $P\in X\cap Y$ if $X\neq Y$
or $P$ a singular point of $X$ if $X=Y$ we decompose $(X)_{P}=(X_{1})_{P}%
\cup...\cup(X_{r})_{P}$, $(Y)_{P}=(Y_{1})_{P}\cup...\cup(Y_{s})_{P}$ into
irreducible curve-germs. Then
\[
C_{P}(X,Y)=\bigcup_{i,j}C_{P}(X_{i},Y_{j}).
\]
Each $C_{P}(X_{i},Y_{j})$ is described in the following way:

\noindent(i) if $(X_{i})_{P}=(Y_{j})_{P}$ and this germ is nonsingular, then
\[
C_{P}(X_{i},Y_{j})=C_{P}(X_{i})=C_{P}(Y_{j})=T_{P}X_{i}=T_{P}Y_{j},
\]

\noindent(ii) if $(X_{i})_{P}\neq(Y_{j})_{P}$ or one of these germs is
singular, then

(1) if $C_{P}(X_{i})\cap C_{P}(Y_{j})=\{P\}$ then
\[
C_{P}(X_{i},Y_{j})=\operatorname*{Span}(C_{P}(X_{i}),C_{P}(Y_{j}))\text{,}%
\]

(2) if $C_{P}(X_{i})=C_{P}(Y_{j})$ then
\begin{align*}
C_{P}(X_{i},Y_{j})  & =\bigcup_{l=1}^{m}\operatorname*{Span}(C_{P}%
(X_{i}),\overline{PQ_{l}})\text{,}\\
1  & \leqslant m\leqslant\min(\deg_{P}X_{i},\deg_{P}Y_{j})
\end{align*}
where $Q_{l}:=\varphi^{-1}(\varphi(P)+\mathbf{v}_{l})$ ($\varphi
:U\rightarrow\mathbb{C}^{n}$ is a canonical map of $\mathbb{P}^{n}$ such that
$P\in U$) and $\mathbf{v}_{l}$ are calculated from local parametrization of
the curves $\varphi(X_{i})-\varphi(P)$ and $\varphi(Y_{j})-\varphi(P)$ at
$\mathbf{0,}$ as it is described in Theorem \ref{opis} (after a linear change
of coordinates in $\mathbb{C}^{n})$.
\end{theorem}

\emph{Acknowledgements. }I thank J.Ch\c{a}dzy\'{n}ski, Z.Jelonek, T. Rodak and
S.Spodzie-ja for helpful comments.

\bigskip

\noindent\textsc{faculty of mathematics}

\noindent\textsc{University of lodz }

\noindent\textsc{ul. s.banacha 22 }

\noindent\textsc{90-238 lodz, Poland}

\bigskip

\noindent\textit{E-mail: }krasinsk$@$krysia.uni.lodz.pl
\end{document}